\documentclass[10pt, A4paper]{article}
\usepackage{amsfonts}
\usepackage{amsmath}
\usepackage{amssymb}

\newtheorem{thm}{Theorem}

\newtheorem{lem}[thm]{Lemma}

\newtheorem{defn}{Definition}
\newtheorem{rem}[thm]{Remark}

\input{tcilatex}
\begin{document}

\title{\textbf{Non-existence of positive weak solutions for some nonlinear }$%
\mathbf{(p,q)}$\textbf{-Laplacian systems}}
\author{\textbf{Salah. A. Khafagy\thanks{\textbf{Current Address}:
Mathematics Department, Faculty of Science in Zulfi, Majmaah University,
Zulfi 11932, P.O. Box 1712, Saudi Arabia.}} \\
Mathematics Department, Faculty of Science, Al-Azhar University,\\
Nasr City (11884), Cairo, Egypt. }
\date{E-mail: el\_gharieb@hotmail.com}
\maketitle

\begin{abstract}
In this work we deal with the class of nonlinear $\mathbf{(p,q)}$\textbf{-}%
Laplacian system of the form%
\begin{equation*}
\left. 
\begin{array}{cc}
-\Delta _{p}u=\mu \rho _{1}(x)f(v) & \text{in}\,\ \Omega , \\ 
-\Delta _{q}v=\nu \rho _{2}(x)g(u) & \text{in}\,\ \Omega , \\ 
u=v=0 & \text{on}\,\ \ \partial \Omega .%
\end{array}%
\right\} 
\end{equation*}

where $\Delta _{p}$ with $p>1$ denotes the $p$-Laplacian defined by $\Delta
_{p}u\equiv div[|\nabla u|^{p-2}\nabla u],$ $\mu ,\nu $ are positive
parameters, $\rho _{1}(x)$, $\rho _{2}(x)$ are weight functions, $%
f,g:[0,\infty )\rightarrow 
\mathbb{R}
$ are continuous functions and $\Omega \subset 
\mathbb{R}
^{N}$ is a bounded domain with smooth boundary $\partial \Omega .$
Non-existence results of positive weak solutions are established under some
certian conditions on $f,g$ when $\mu \nu $ is small.
\end{abstract}

\date{\textbf{2000 Mathematics Subject Classification}: 35J60,35B40.\\
\textbf{Key words}: weak solution, p-Laplacian.}

\section{Introduction:}

In this paper we first consider a non-existence result of positive weak
solutions for the following nonlinear system%
\begin{equation}
\left. 
\begin{array}{cc}
-\Delta _{p}u=\lambda a_{1}(x)v^{p-1}-b_{1}(x)v^{\alpha -1}-c_{1}(x) & \text{%
in}\,\ \Omega , \\ 
-\Delta _{q}v=\lambda a_{2}(x)u^{q-1}-b_{2}(x)u^{\beta -1}-c_{2}(x) & \text{%
in}\,\ \Omega , \\ 
u=v=0 & \text{on}\,\ \ \partial \Omega ,%
\end{array}%
\right\}   \tag{1.1}  \label{11}
\end{equation}%
where $\Delta _{p}$ with $p>1$ denotes the weighted $p$-Laplacian defined by 
$\Delta _{p}u\equiv div[|\nabla u|^{p-2}\nabla u],$ $a_{i}(x)$, $b_{i}(x)$
and $c_{i}(x),i=1,2$ are weight functions, $\alpha $ and $\beta $ are
positive constants and $\Omega \subset 
\mathbb{R}
^{N}$ is a bounded domain with smooth boundary $\partial \Omega .$

We first show that if $\lambda <\max (\lambda _{p},\lambda _{q})$, where $%
\lambda _{p},\lambda _{q}$ is the first eigenvalue of $-\Delta _{p},-\Delta
_{q}$ respectively, then system (\ref{11}) has no positive weak solutions.

Next we consider the nonlinear system%
\begin{equation}
\left. 
\begin{array}{cc}
-\Delta _{p}u=\mu \rho _{1}(x)f(v) & \text{in}\,\ \Omega , \\ 
-\Delta _{q}v=\nu \rho _{2}(x)g(u) & \text{in}\,\ \Omega , \\ 
u=v=0 & \text{on}\,\ \ \partial \Omega .%
\end{array}%
\right\}  \tag{1.2}  \label{12}
\end{equation}%
where $\mu ,\nu $ are positive parameters, $\rho _{1}(x)$, $\rho _{2}(x)$
are weight functions and $\Omega \subset 
\mathbb{R}
^{N}$ is a bounded domain with smooth boundary $\partial \Omega $. Let $%
f,g:[0,\infty )\rightarrow 
\mathbb{R}
$ are continuous functions. Also, assume that there exist positive numbers $%
K_{i}$ and $M_{i},i=1,2$ such that%
\begin{equation}
f(v)\leq K_{1}v^{p-1}-M_{1},\text{ \ \ \ for all\ \ }v\geq 0  \tag{1.3}
\label{13}
\end{equation}%
and%
\begin{equation}
g(u)\leq K_{2}u^{q-1}-M_{2},\text{ \ \ for all\ \ }u\geq 0.  \tag{1.4}
\label{14}
\end{equation}

We discuss a non-existence result for system (\ref{12}) when $\mu \nu $ is
small.

Problems of the form (1.1) and (\ref{12}) arise from many branches of pure
mathematics as in the theory of quasiregular and quasiconformal mappings
(see \cite{T1984}) as well as from various problems in mathematical physics
notably the flow of non-Newtonian fluids. The $p$-Laplacian also appears in
the study of torsional creep (elastic for $p=2$, plastic as $p\rightarrow
\infty $, (see \cite{K1990}), glacial sliding ($p\in (1;\frac{4}{3}]$, see 
\cite{PR1974} or flow through porous media ($p=\frac{3}{2}$, see \cite%
{SW1991}). For existence and non-existence results of positive weak
solutions for systems involving the weighted $p$-Laplacian, see (\cite%
{C2001,DH2001,K2011,K2012,K2013,K2009,SK2009}).

This paper is organized as follows: In section 2, we introduce some
technical results and notations, which are established in \cite{DKN1997}. In
section 3, we prove the non-existence of positive weak solutions for system (%
\ref{11}) and (\ref{12}).

\section{Technical Results}

Let us introduce the Sobolev space $W^{1,p}(\Omega ),$ $1<p<\infty ,$
defined as the completion of $C^{\infty }(\Omega )$ with respect to the norm
(see \cite{DKN1997}) 
\begin{equation}
\Vert u\Vert _{W^{1,p}(\Omega )}=\left[ \int\limits_{\Omega
}|u|^{p}+\int\limits_{\Omega }|\nabla u|^{p}\right] ^{\frac{1}{p}}<\infty . 
\tag{2.1}  \label{21}
\end{equation}

Since we are dealing with the Dirichlet problem, we define the space $%
W_{0}^{1,p}(\Omega )$ as the closure of $C_{0}^{\infty }(\Omega )$ in $%
W^{1,p}(\Omega )$ with respect to the norm%
\begin{equation}
\Vert u\Vert _{W_{0}^{1,p}(\Omega )}=\left[ \int\limits_{\Omega }|\nabla
u|^{p}\right] ^{\frac{1}{p}}<\infty ,  \tag{2.2}  \label{22}
\end{equation}%
which is equivalent to the norm given by (\ref{21}). Both spaces $%
W^{1,p}(\Omega )$ and $W_{0}^{1,p}(\Omega )$ are well defined reflexive
Banach Spaces.

Now, we introduce some technical results concerning the eigenvalue problem 
\begin{equation}
\left. 
\begin{array}{cc}
-\Delta _{p}u=\lambda a(x)|u|^{p-2}u & \text{in}\,\ \Omega , \\ 
u=0 & \text{on }\,\partial \Omega .%
\end{array}%
\right\}  \tag{2.3}  \label{23}
\end{equation}

We will say $\lambda \in R$ is an eigenvalue of (\ref{23}) if there exists $%
u\in W_{0}^{1,p}(\Omega ),$ $u\neq 0,$ such that

\begin{equation}
\int\limits_{\Omega }|\nabla u|^{p-2}\nabla u\nabla \phi dx=\lambda
\int\limits_{\Omega }a(x)u^{p-2}u\phi dx,  \tag{2.4}  \label{24}
\end{equation}

holds for $\phi \in W_{0}^{1,p}(\Omega ).$ Then u is called an eigenfunction
corresponding to the eigenvalue $\lambda .$

\begin{lem}
There exists the first eigenvalue $\lambda _{p}>0$ and precisely one
corresponding eigenfunction $\phi _{p}\geq 0\ $a.e.\ in$\ \Omega $ of the
eigenvalue problem (\ref{23}). Moreover, it is characterized by 
\begin{equation*}
\lambda _{p}=\frac{\int\limits_{\Omega }|\nabla \phi _{p}|^{p}}{%
\int\limits_{\Omega }a(x)|\phi _{p}|^{p}}=\inf_{u\in W_{0}^{1,p}(\Omega )}%
\frac{\int\limits_{\Omega }|\nabla u|^{p}}{\int\limits_{\Omega }a(x)|u|^{p}}%
\leq \frac{\int\limits_{\Omega }|\nabla u|^{p}}{\int\limits_{\Omega
}a(x)|u|^{p}}=\lambda .
\end{equation*}
\end{lem}

\begin{defn}
A pair of non-negative functions $(u,v)\in W_{0}^{1,p}(\Omega )\times
W_{0}^{1,q}(\Omega )$ are called a weak solution of (\ref{12}) if they
satisfy 
\begin{eqnarray*}
\int\limits_{\Omega }|\nabla u|^{p-2}\nabla u\nabla \zeta dx &=&\mu
\int\limits_{\Omega }\rho _{1}(x)f(v)\zeta dx, \\
\int\limits_{\Omega }|\nabla v|^{q-2}\nabla v\nabla \eta dx &=&\nu
\int\limits_{\Omega }\rho _{2}(x)g(u)\eta dx,
\end{eqnarray*}%
for all test functions $\zeta \in W_{0}^{1,p}(P,\Omega ),\eta \in
W_{0}^{1,q}(\Omega ).$
\end{defn}

\section{Non-existence Results}

In this section we state our main results. Throught this section, we assume $%
q$ be such that $\frac{1}{p}+\frac{1}{q}=1.$

\begin{thm}
For $\lambda \leq \lambda ^{\ast }$, system (\ref{11}) has no positive weak
solution.
\end{thm}

\textbf{Proof. }Assume that there exist a positive solution $(u,v)\in
W_{0}^{1,p}(\Omega )\times W_{0}^{1,q}(\Omega )$ of (\ref{11}). Multiplying
the first equation of ( \ref{11}) by $u,$ we have%
\begin{eqnarray}
\int\limits_{\Omega }|\nabla u|^{p}dx &=&\int\limits_{\Omega }[\lambda
a_{1}(x)v^{p-1}-b_{1}(x)v^{\alpha -1}-c_{1}(x)]udx  \TCItag{3.1}  \label{31}
\\
&<&\int\limits_{\Omega }[\lambda a_{1}(x)v^{p-1}-c_{1}(x)]udx.  \notag
\end{eqnarray}

But, from the characterization of the first eigenvalue, we have%
\begin{equation}
\lambda _{p}\int\limits_{\Omega }a(x)|u|^{p}\leq \int\limits_{\Omega
}|\nabla u|^{p}.  \tag{3.2}  \label{32}
\end{equation}

Combining (\ref{31}) and (\ref{32}), we have%
\begin{equation}
\lambda _{p}\int\limits_{\Omega }a_{1}(x)u^{p}<\int\limits_{\Omega }\lambda
a_{1}(x)v^{p-1}udx-\int\limits_{\Omega }c_{1}(x)udx.  \tag{3.3}  \label{33}
\end{equation}

Similarly, from the second equation of \ (\ref{11}), we obtain%
\begin{equation}
\lambda _{q}\int\limits_{\Omega }a_{2}(x)v^{q}<\int\limits_{\Omega }\lambda
a_{2}(x)u^{q-1}vdx-\int\limits_{\Omega }c_{2}(x)vdx.  \tag{3.4}  \label{34}
\end{equation}

Adding (\ref{33}) and (\ref{34}), we get%
\begin{eqnarray*}
\lambda _{p}\int\limits_{\Omega }a_{1}(x)u^{p}+\lambda
_{q}\int\limits_{\Omega }a_{2}(x)v^{q} &<&\int\limits_{\Omega }\lambda
a_{1}(x)v^{p-1}udx+\int\limits_{\Omega }\lambda a_{2}(x)u^{q-1}vdx \\
&&-\int\limits_{\Omega }c_{1}(x)udx-\int\limits_{\Omega }c_{2}(x)vdx \\
&<&\int\limits_{\Omega }\lambda a_{1}(x)v^{p-1}udx+\int\limits_{\Omega
}\lambda a_{2}(x)u^{q-1}vdx
\end{eqnarray*}

Applying the Young inequality on the right hand side of the above equation,
we have%
\begin{equation}
\lambda _{p}\int\limits_{\Omega }a_{1}(x)u^{p}+\lambda
_{q}\int\limits_{\Omega }a_{2}(x)v^{q}<\int\limits_{\Omega }\lambda a_{1}(x)[%
\frac{u^{p}}{p}+\frac{v^{p}}{q}]dx+\int\limits_{\Omega }\lambda a_{2}(x)[%
\frac{v^{q}}{q}+\frac{u^{q}}{p}]dx  \tag{3.5}  \label{35}
\end{equation}

Now, we discuss the following two cases:

Case I, if $u\leq v$ for all $x,$ then (\ref{35}) becomes

\begin{equation*}
\lambda _{p}\int\limits_{\Omega }a_{1}(x)v^{p}+\lambda
_{q}\int\limits_{\Omega }a_{2}(x)v^{q}<\int\limits_{\Omega }\lambda
a_{1}(x)v^{p}dx+\int\limits_{\Omega }\lambda a_{2}(x)v^{q}dx.
\end{equation*}

Hence,

\begin{equation*}
(\lambda _{p}-\lambda )\int\limits_{\Omega }a_{1}(x)v^{p}+(\lambda
_{q}-\lambda )\int\limits_{\Omega }a_{2}(x)v^{q}<0
\end{equation*}

which is a contradiction if $\lambda \leq \max (\lambda _{p},\lambda
_{q})=\lambda ^{\ast }.$

Case II, if $u\geq v$ \ for all $x,$ then (\ref{35}) becomes

\begin{equation*}
\lambda _{p}\int\limits_{\Omega }a_{1}(x)u^{p}+\lambda
_{q}\int\limits_{\Omega }a_{2}(x)u^{q}<\int\limits_{\Omega }\lambda
a_{1}(x)u^{p}dx+\int\limits_{\Omega }\lambda a_{2}(x)u^{q}dx.
\end{equation*}

Hence,

\begin{equation*}
(\lambda _{p}-\lambda )\int\limits_{\Omega }a_{1}(x)u^{p}+(\lambda
_{q}-\lambda )\int\limits_{\Omega }a_{2}(x)u^{q}<0
\end{equation*}

which is a contradiction if $\lambda \leq \max (\lambda _{p},\lambda
_{q})=\lambda ^{\ast }.$ The proof complete.

Now we consider the main result for system (\ref{22}):

\begin{thm}
Let (\ref{13}) and (\ref{14}) hold. Then system (\ref{12}) has no positive
weak solution if $\mu \nu \leq \frac{\lambda _{1}^{2}}{K_{1}K_{2}}.$
\end{thm}

\textbf{Proof.} Suppose $u>0$ and $v>0$ be such that $(u,v)$ is a solution
of (\ref{22}). We prove our theorem by arriving at a contradiction.
Multiplying the first equation in (\ref{22}) by a positive eigenfunction say 
$\phi _{p}$ corresponding to $\lambda _{p}$, we obtain%
\begin{equation*}
-\int\limits_{\Omega }\Delta _{p}u\text{ }\phi _{p}dx=\mu
\int\limits_{\Omega }\rho _{1}(x)f(v)\phi _{p}dx,
\end{equation*}

and hence using (\ref{21}) and (\ref{13}), we have%
\begin{equation}
\lambda _{p}\int\limits_{\Omega }\rho _{1}(x)u^{p-1}\phi _{p}dx\leq \mu
\int\limits_{\Omega }\rho _{1}(x)[K_{1}v^{q-1}-M_{1}]\phi _{p}dx.  \tag{3.6}
\label{36}
\end{equation}

Similarly using the second equation in (\ref{22}) and (\ref{14}) we obtain%
\begin{equation}
\lambda _{q}\int\limits_{\Omega }\rho _{2}(x)v^{q-1}\phi _{q}dx\leq \nu
\int\limits_{\Omega }\rho _{2}(x)[K_{2}u^{p-1}-M_{2}]\phi _{q}dx.  \tag{3.7}
\label{37}
\end{equation}

From (\ref{37}), we have%
\begin{equation}
v^{q-1}\leq \frac{\nu }{\lambda _{q}}[K_{2}u^{p-1}-M_{2}]  \tag{3.8}
\label{38}
\end{equation}

Combining (\ref{36}) and (\ref{38}) we obtain%
\begin{equation*}
\lbrack \lambda _{p}-\mu \nu \frac{K_{1}K_{2}}{\lambda _{q}}%
]\int\limits_{\Omega }\rho _{1}(x)u^{p-1}\phi _{p}\leq -\mu
\int\limits_{\Omega }\rho _{1}(x)[\frac{\nu K_{1}M_{2}}{\lambda _{q}}%
+M_{1}]\phi _{p}<0.
\end{equation*}

Hence system (\ref{12}) has no positive weak solution if $\mu \nu \leq \frac{%
\lambda _{p}\lambda _{q}}{K_{1}K_{2}}.$

\begin{rem}
If $f,g$ be such that

\begin{equation}
f(v)\geq K_{1}v^{q-1}+M_{1},\text{ \ \ \ for all\ \ }v\geq 0,  \tag{3.9}
\label{39}
\end{equation}%
and%
\begin{equation}
g(u)\geq K_{2}u^{p-1}+M_{2},\text{ \ \ for all\ \ }u\geq 0,  \tag{3.10}
\label{310}
\end{equation}%
then we have the following theorem:

\begin{thm}
Let (\ref{39}) and (\ref{310}) hold. Then system (\ref{12}) has no positive
weak solution if $\mu \nu \geq \frac{\lambda _{p}\lambda _{q}}{K_{1}K_{2}}.$
\end{thm}
\end{rem}

\textbf{Proof. }The proof proceeds in the same way as for Theorem 6.

\begin{rem}
When $p=q$, $m_{1}(x)=m_{2}(x)=m,m=a,b,c$ is constant , and $\alpha =\beta ,$
we have some results for (\ref{11}) in \cite{A2011}.
\end{rem}

\begin{rem}
When $p=q$ and\ $\rho _{1}(x)=\rho _{2}(x)=1,$ we have some results for (\ref%
{12}) in \cite{A2011} .
\end{rem}

$\mathbf{Acknowledgement.}$ The author would like to express his gratitude
to Professor H. M. Serag (Mathematics Department, Faculty of Science, AL-
Azhar University) for continuous encouragement during the development of
this work.

\end{document}